\documentclass[preprint,3p,times]{elsarticle}

\usepackage{amssymb}
\usepackage{amsmath}
\usepackage[all,cmtip]{xy}
\usepackage{amsthm}
\usepackage{color}
\usepackage{enumitem}
\usepackage{tikz}
\usepackage{amsmath, mathrsfs}
\usepackage{appendix}
\usepackage[colorlinks=true]{hyperref}
\usepackage[hyphenbreaks]{breakurl}
\input diagxy

\theoremstyle{plain}
  \newtheorem{thm}{Theorem}[section]
  \newtheorem{lem}[thm]{Lemma}
  \newtheorem{prop}[thm]{Proposition}
  \newtheorem{cor}[thm]{Corollary}
\theoremstyle{definition}
  \newtheorem{defn}[thm]{Definition}
  \newtheorem{exmp}[thm]{Example}
  \newtheorem{rem}[thm]{Remark}
  
  \newtheorem{ques}[thm]{Question}

\DeclareMathAlphabet{\mathcal}{OMS}{cmsy}{m}{n}

\DeclareMathOperator{\id}{id}

\makeatletter
\def\ps@pprintTitle{%
 \let\@oddhead\@empty
 \let\@evenhead\@empty
 \def\@oddfoot{\centerline{\thepage}}%
 \let\@evenfoot\@oddfoot}
\makeatother

\def\oto{{\bfig\morphism<180,0>[\mkern-4mu`\mkern-4mu;]\place(86,0)[\circ]\efig}}
\def\rto{{\bfig\morphism<180,0>[\mkern-4mu`\mkern-4mu;]\place(82,0)[\mapstochar]\efig}}

\newcommand{\ra}{\rightarrow}

\newcommand{\lda}{\swarrow}
\newcommand{\rda}{\searrow}

\newcommand{\bv}{\bigvee}
\newcommand{\bw}{\bigwedge}

\newcommand{\dv}{\dashv}

\newcommand{\nat}{\natural}

\newcommand{\dbv}{\displaystyle\bv}

\renewcommand{\phi}{\varphi}

\newcommand{\ep}{\varepsilon}

\newcommand{\Si}{\Sigma}

\newcommand{\CC}{\mathcal{C}}

\newcommand{\CQ}{\mathcal{Q}}

\newcommand{\sF}{\mathsf{F}}

\newcommand{\sK}{\mathsf{K}}

\newcommand{\sP}{\mathsf{P}}
\newcommand{\sQ}{\mathsf{Q}}

\newcommand{\sT}{\mathsf{T}}
\newcommand{\sU}{\mathsf{U}}

\newcommand{\sm}{\mathsf{m}}

\newcommand{\sfs}{\mathsf{s}}

\newcommand{\bbZ}{\mathbb{Z}}

\newcommand{\Map}{\mathbf{Map}}

\newcommand{\Rel}{\mathbf{Rel}}
\newcommand{\Set}{\mathbf{Set}}

\newcommand{\QMap}{\sQ\text{-} \mathbf{Map}}
\newcommand{\QParMap}{\sQ\text{-}\Map^{\partial}}

\newcommand{\sQMap}{\{\star\}/\QMap}

\newcommand{\QRel}{\sQ\text{-}\Rel}

\newcommand{\op}{{\rm op}}

\newcommand{\PM}{\sP M}
\newcommand{\PX}{\sP X}

\newcommand{\with}{\mathrel{\&}}

\renewcommand{\leq}{\leqslant}
\renewcommand{\geq}{\geqslant}

\numberwithin{equation}{section}

\allowdisplaybreaks

\begin{document}

\begin{frontmatter}



\title{Quantale-valued maps and partial maps}


\author{Lili Shen\corref{cor}}
\ead{shenlili@scu.edu.cn}

\author{Xiaoye Tang}
\ead{tang.xiaoye@qq.com}

\cortext[cor]{Corresponding author.}
\address{School of Mathematics, Sichuan University, Chengdu 610064, China}

\begin{abstract}
Let $\mathsf{Q}$ be a commutative and unital quantale. By a $\mathsf{Q}$-map we mean a left adjoint in the quantaloid of sets and $\mathsf{Q}$-relations, and by a partial $\mathsf{Q}$-map we refer to a Kleisli morphism with respect to the maybe monad on the category $\mathsf{Q}\text{-}\mathbf{Map}$ of sets and $\mathsf{Q}$-maps. It is shown that every $\mathsf{Q}$-map is symmetric if and only if $\mathsf{Q}$ is weakly lean, and that every $\mathsf{Q}$-map is exactly a map in $\mathbf{Set}$ if and only $\mathsf{Q}$ is lean. Moreover, assuming the axiom of choice, it is shown that the category of sets and partial $\mathsf{Q}$-maps is monadic over $\mathsf{Q}\text{-}\mathbf{Map}$.
\end{abstract}


\end{frontmatter}




\section{Introduction}

Given (crisp) sets $X$ and $Y$, a \emph{partial map} $f$ from $X$ to $Y$ is a map from a (possibly empty) subset $X'\subseteq X$ to $Y$. For the category
\[\Set^{\partial}\] 
of (crisp) sets and partial maps, the following results are well known \cite{MacLane1998,Riehl2016}:
\begin{itemize}
\item $\Set^{\partial}$ is equivalent to the coslice category $\{\star\}/\Set$, where $\{\star\}$ is the singleton set.
\item The forgetful functor $\sU\colon\{\star\}/\Set\to\Set$ admits a left adjoint, which carries a set $X$ to the inclusion map
\[\{\star\}\ \to/^(->/ X\amalg\{\star\},\]
where $X\amalg\{\star\}$ refers to the disjoint union of $X$ and the singleton set.
\item The induced monad on $\Set$ is called the \emph{maybe monad}, whose Eilenberg-Moore and Kleisli categories are $\{\star\}/\Set$ and $\Set^{\partial}$, respectively (cf. \cite[Examples 5.1.4 and 5.2.10]{Riehl2016}). 
\end{itemize}
Therefore, the Eilenberg-Moore and Kleisli categories of the maybe monad on $\Set$ are equivalent. In particular, $\{\star\}/\Set$ and $\Set^{\partial}$ are both monadic over $\Set$.

It is then natural to ask whether it is possible to formulate the notion of partial map, and establish the above results in the fuzzy setting. To achieve this goal, we note that there is a very general framework in category theory for partial maps; see, e.g., \cite[Exercise IV.2.11]{MacLane1998}. Explicitly, let $\CC$ be a category with finite coproducts. For every $\CC$-object $A$, the forgetful functor from the coslice category $A/\CC$ to $\CC$, i.e.,
\begin{equation} \label{U-AC-C}
\sU\colon A/\CC\to\CC,\quad (A\ra C)\mapsto C,
\end{equation}
admits a left adjoint, given by
\begin{equation} \label{A-AC}
C\mapsto(A\ra A\amalg C).
\end{equation}
It is not difficult to see that the Eilenberg-Moore category of the induced monad is isomorphic to $A/\CC$, and thus $A/\CC$ is monadic over $\CC$. The maybe monad on $\Set$ is simply obtained by setting $\CC=\Set$ and $A=\{\star\}$.

Considering a non-trivial, commutative and unital \emph{quantale}  \cite{Rosenthal1990} $\sQ$ as the table of truth values and following the terminology of \cite[Definition 2.3.1]{Heymans2010}, by a \emph{$\sQ$-map} (i.e., a map valued in the quantale $\sQ$) we mean a left adjoint in the \emph{quantaloid} \cite{Rosenthal1996,Stubbe2014}
\[\QRel\]
of sets and $\sQ$-relations; that is, a $\sQ$-relation $\zeta\colon X\rto Y$ admitting a right adjoint in $\QRel$, with the value
\[\zeta(x,y)\]
interpreted as the extent of $y$ being the image of $x$ under the map $\zeta$ (see Remark \ref{why-Q-map} for detailed elaboration). In Section \ref{Quantale-valued_maps} we prove the following results (Theorems \ref{Q-map-symmetric-weakly-lean} and \ref{QMap-Set-lean}), which reveal that for a general quantale $\sQ$, $\sQ$-maps are essentially different from crisp maps between sets:
\begin{itemize}
\item Every $\sQ$-map is \emph{symmetric} (in the sense that the right adjoint is given by its opposite) if, and only if, $\sQ$ is \emph{weakly lean}. This is a discrete counterpart of \cite[Proposition 3.5.3]{Heymans2010} and \cite[Theorem 3.7]{Heymans2011}.
\item Every $\sQ$-map is the graph of a map in $\Set$ (so that the category $\QMap$ of sets and $\sQ$-maps is isomorphic to $\Set$) if, and only if, $\sQ$ is \emph{lean}. This is a generalization of \cite[Proposition III.1.2.1]{Hofmann2014} to the non-integral setting.
\end{itemize}

With the above preparations we are now able to postulate the notion of \emph{partial $\sQ$-map} through the \emph{maybe monad} on the category $\QMap$, which is obtained by setting $\CC=\QMap$ and $A=\{\star\}$ in \eqref{U-AC-C} and \eqref{A-AC}. Explicitly, by a partial $\sQ$-map $\zeta$ from $X$ to $Y$ we refer to a $\sQ$-map
\[\zeta\colon X\rto Y\amalg\{\star\},\]
which is precisely a morphism in the Kleisli category of the maybe monad on $\QMap$, with the value 
\[\zeta(x,\star)\]
interpreted as the degree that $x$ has no image under the map $\zeta$ (see Remark \ref{why-partial-Q-map} for detailed elaboration). The main result of Section \ref{Quantale-valued_partial_maps}, Theorem \ref{QMap-iso-sQMap}, is to show that assuming the axiom of choice, the Eilenberg-Moore category (i.e., the coslice category $\{\star\}/\QMap$) and the Kleisli category (i.e., the category $\QParMap$ of sets and partial $\sQ$-maps) of the maybe monad on $\QMap$ are equivalent, so that $\sQMap$ and $\QParMap$ are both monadic over $\QMap$. It is noteworthy to point out that this result may not hold in the general categorical setting (i.e., the monad induced by \eqref{U-AC-C} and \eqref{A-AC} for a general category $\CC$), and an equivalent statement of the axiom of choice given in \cite{Hajnal1972} (i.e., every non-empty set admits a group structure) plays a vital role in its proof.

\section{Quantales}

Throughout, let 
\[\sQ=(\sQ,\with,k)\]
denote a non-trivial, commutative and unital \emph{quantale} \cite{Mulvey1986,Rosenthal1990,Eklund2018}, which consists of a complete lattice $\sQ$ (with the top and bottom elements denoted by $\top$ and $\bot$, respectively), an element $k\in\sQ$ and a binary operation $\with$ on $\sQ$, such that 
\begin{itemize}
\item $\bot<k$,
\item $(\sQ,\with,k)$ is a commutative monoid with unit $k$, and
\item $p\with \Big(\dbv\limits_{i\in I}q_i\Big)=\dbv\limits_{i\in I}p\with q_i$
for all $p, q_{i}\in \sQ\,(i\in I)$.
\end{itemize}

The right adjoint induced by the multiplication $\with$, denoted by $\ra$,
\[(p\with -)\dashv(p\ra -)\colon\sQ\to\sQ,\]
satisfies
\[p\with q\leq r\iff p\leq q\ra r\]
for all $p,q,r\in\sQ$. 

We say that $\sQ$ is \emph{integral} if $k=\top$, and $\sQ$ is \emph{divisible} \cite{Hoehle1995a,Hajek1998} if 
\[d\leq q\implies \exists p\in\sQ\colon d=p\with q\]
for all $d,q\in\sQ$. In particular, divisible quantales are necessarily integral (see, e.g., \cite[Proposition 2.1]{Pu2012}).

\begin{exmp} \label{quantale-example}
We list here some examples of commutative and unital quantales that will be concerned later:
\begin{enumerate}[label=(\arabic*)]
\item \label{quantale-example:C3} On the three-chain $C_3=\{\bot,k,\top\}$ we have the quantale $(C_3,\with,k)$, with 
\[\top\with\top=\top\]
and the other multiplications being trivial.
\item \label{quantale-example:frame} Every frame is a divisible quantale. In particular: 
\begin{itemize}
\item The partially ordered sets 
\[F_1=\{\bot,p,q,\top\}\quad\text{and}\quad F_2=\{\bot,p,q,r,\top\}\]
illustrated by the Hasse diagrams
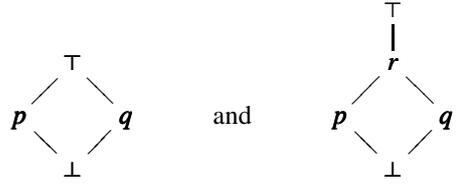
\[\bfig
\Atriangle/-`-`/<200,200>[\top`p`q;``]
\Vtriangle(0,-200)/`-`-/<200,200>[p`q`\bot;``]
\place(800,0)[\text{and}]
\Atriangle(1200,0)/-`-`/<200,200>[r`p`q;``]
\morphism(1400,200)/-/<0,200>[r`\top;]
\Vtriangle(1200,-200)/`-`-/<200,200>[p`q`\bot;``]
\place(1700,0)[\text{,}]
\efig\]
respectively, are both frames. 
\item For each (crisp) set $X$, the powerset $\PX=\{A\mid A\subseteq X\}$ equipped with the inclusion order is a frame.
\end{itemize}
\item \label{quantale-example:diamond}
On the diamond lattice $M_3$ given by the Hasse diagram
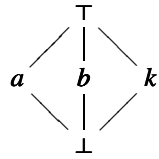
\[\bfig
\Atrianglepair/-`-`-``/<250,250>[\top`a`b`k;````]
\Vtrianglepair(0,-250)/``-`-`-/<250,250>[a`b`k`\bot;````]
\efig\]
we have the quantales $M_3=(M_3,\with,k)$ and $M'_3=(M_3,\with',k)$ with
\[a\with a=b,\quad b\with b=a,\quad a\with b=k,\quad a\with\top=b\with\top=\top\]
and
\[a\with' a=a\with' b=b\with' b=a\with'\top=b\with'\top=\top.\]
\item \label{quantale-example:t-norm} Every continuous t-norm \cite{Klement2000,Klement2004b,Alsina2006} $*$ on the unital interval $[0,1]$ gives rise to a divisible quantale $([0,1],*,1)$. 
\item \label{quantale-example:Lawvere} Let $[-\infty,\infty]$ be the extended real line equipped with the order ``$\geq$''. Then $([-\infty,\infty],+,0)$ is a quantale, which includes the Lawvere quantale \cite{Lawvere1973} $([0,\infty],+,0)$ as a divisible subquantale.
\item \label{quantale-example:free} Every commutative monoid $(M,\with,k)$ induces a \emph{free quantale} $(\PM,\with,\{k\})$, where $\PM$ is the powerset of $M$, and
\[A\with B=\{a\with b\mid a\in A,\ b\in B\}\]
for all $A,B\subseteq M$.
\end{enumerate}
\end{exmp}

A \emph{quantaloid} \cite{Rosenthal1996,Stubbe2005,Stubbe2014} $\CQ$ is a category with a class of objects $\CQ_0$ in which $\CQ(p,q)$ is a complete lattice for all $p,q\in\CQ_0$, such that the composition $\circ$ of morphisms preserves joins on both sides, i.e.,
\[v\circ\Big(\bv\limits_{i\in I}u_i\Big)=\bv\limits_{i\in I}(v\circ u_i)\quad\text{and}\quad \Big(\bv\limits_{i\in I} v_i\Big)\circ u=\bv\limits_{i\in I}(v_i\circ u)\]
for all $u,u_i\in\CQ(p,q)$, $v,v_i\in\CQ(q,r)$ $(i\in I)$. The corresponding right adjoints induced by the compositions 
\[(-\circ u)\dv(-\lda u)\colon\CQ(p,r)\to\CQ(q,r)\quad\text{and}\quad(v\circ -)\dv(v\rda -)\colon\CQ(p,r)\to\CQ(p,q)\]
satisfy 
\[v\circ u\leq w\iff v\leq w\lda u\iff u\leq v\rda w\] 
for all $\CQ$-arrows $u\colon p\to q$, $v\colon q\to r$, $w\colon p\to r$.

An \emph{adjunction} in a quantaloid $\CQ$ is a pair of $\CQ$-arrows $u\colon p\to q$ and $v\colon q\to p$, denoted by $u\dv v$, such that
\[1_p\leq v\circ u\quad\text{and}\quad u\circ v\leq 1_q.\]

\begin{prop} 
If $u\dv v\colon q\to p$ in a quantaloid $\CQ$, then
\[v=u\rda 1_q\quad\text{and}\quad u=1_q\lda v.\]
\end{prop}

\begin{proof}
Note that $1_p\leq v\circ u$ implies that $u\rda 1_q\leq v\circ u\circ(u\rda 1_q)\leq v$, and the reverse inequality is an immediate consequence of $u\circ v\leq 1_q$. Therefore $v=u\rda 1_q$, and $u=1_q\lda v$ can be proved analogously.
\end{proof}

Therefore, the right adjoint of a $\CQ$-arrow, when it exists, is necessarily unique. In what follows we define 
\begin{equation} \label{u*-def}
u^*:=u\rda 1_q\colon q\to p
\end{equation}
for each $\CQ$-arrow $u\colon p\to q$. $u$ is called a \emph{map} in $\CQ$ \cite{Heymans2010} if $u\dv u^*$, and we denote by
\[\Map(\CQ)\]
the subcategory of $\CQ$ whose objects are the same as $\CQ$, and whose morphisms are maps in $\CQ$.

Since it always holds that $u\circ u^*=u\circ(u\rda 1_q)\leq 1_q$, we immediately obtain the following lemma:

\begin{lem} \label{u-map-leq}
A $\CQ$-arrow $u\colon p\to q$ is a map in $\CQ$ if, and only if, $1_p\leq u^*\circ u$.
\end{lem}

Given sets $X$, $Y$, a \emph{$\sQ$-relation} $\phi\colon X\rto Y$ is a function $\phi\colon X\times Y\to\sQ$. We denote by
\[\phi^{\op}:Y\rto X,\quad \phi^{\op}(y,x)=\phi(x,y)\]
the \emph{opposite} of a $\sQ$-relation $\phi\colon X\rto Y$. Sets and $\sQ$-relations constitute a quantaloid
\[\QRel,\]
in which the local order is inherited from $\sQ$, and 
\begin{align}
& \psi\circ\phi\colon X\rto Z,\quad(\psi\circ\phi)(x,z)=\bv\limits_{y\in Y}\psi(y,z)\with\psi(x,y),\label{QRel-comp}\\
& \xi\lda\phi\colon Y\rto Z,\quad (\xi\lda\phi)(y,z)=\bw\limits_{x\in X}(\phi(x,y)\ra\xi(x,z)),\label{QRel-left-imp}\\
& \psi\rda\xi\colon X\rto Y,\quad (\psi\rda\xi)(x,y)=\bw\limits_{z\in Z}\psi(y,z)\ra\xi(x,z) \label{QRel-right-imp}
\end{align}
for all $\sQ$-relations $\phi\colon X\rto Y$, $\psi\colon Y\rto Z$, $\xi\colon X\rto Z$; the identity $\sQ$-relation on $X$ is given by
\[\id\colon X\rto X,\quad \id_X(x,y)=\begin{cases}
k & \text{if}\ x=y,\\
\bot & \text{else}.
\end{cases}\]
Note that $\sQ$-relations are also known as \emph{$\sQ$-matrices} \cite{Heymans2010}. In fact, a $\sQ$-relation $\phi\colon X\rto Y$ between finite sets may be written as a matrix $(\phi_{yx})_{Y\times X}$, where $\phi_{yx}=\phi(x,y)$. 
Then its composite with $(\psi_{zy})_{Z\times Y}$ can be regarded as the composite of matrices, i.e.,
\begin{equation} \label{matrices}
((\psi\circ\phi)_{zx})_{Z\times X}\quad\text{with}\quad (\psi\circ\phi)_{zx}=\bv\limits_{y\in Y}\psi_{zy}\with\phi_{yx}.
\end{equation}
Given $\sQ$-relations $\phi\colon X\rto Y$ and $\psi\colon Z\rto W$, their disjoint union 
\[\psi\amalg\phi\colon X\amalg Z\rto Y\amalg W\]
is given by
\[(\psi\amalg\phi)(x,y)=\phi(x,y),\quad(\psi\amalg\phi)(z,w)=\psi(z,w)\quad\text{and}\quad(\psi\amalg\phi)(x,w)=(\psi\amalg\phi)(z,y)=\bot\]
for all $x\in X$, $z\in Z$, $y\in Y$, $w\in W$. It is straightforward to verify that
\begin{equation} \label{psi-circ-phi-amalg-xi}
(\psi\circ\phi)\amalg\xi=(\psi\amalg\xi)\circ(\phi\amalg\xi)
\end{equation}
for all $\sQ$-relations $\phi\colon X\rto Y$, $\psi\colon Y\rto Z$, $\xi\colon Z\rto W$.


\section{Quantale-valued maps} \label{Quantale-valued_maps}

\begin{defn} \label{Q-map}
Let $X$, $Y$ be (crisp) sets. A \emph{$\sQ$-map} $\zeta$ from $X$ to $Y$ is a map $\zeta\colon X\rto Y$ in the quantaloid $\QRel$.
\end{defn}

Sets and $\sQ$-maps constitute a category 
\[\QMap:=\Map(\QRel).\]

\begin{rem} \label{why-Q-map}
To justify Definition \ref{Q-map}, let $\zeta\colon X\rto Y$ be a $\sQ$-map. First of all, the value
\[\zeta(x,y)\]
is interpreted as the extent of $y$ being the image of $x$ under the map $\zeta$. By Equations \eqref{u*-def} and \eqref{QRel-right-imp}, the value
\begin{equation} \label{zeta*-pointwise}
\zeta^*(y,x)=(\zeta\rda\id_Y)(y,x)=\bw_{z\in Y}\zeta(x,z)\ra\id_Y(y,z)
\end{equation}
also represents the extent of $y$ being the image of $x$ under the map $\zeta$, since the last expression of \eqref{zeta*-pointwise} may be understood as:
\begin{itemize}
\item For each $z\in Y$, if $z$ is the image of $x$ under $\zeta$, then $z$ is equal to $y$.
\end{itemize}  
Therefore, the adjunction $\zeta\dv\zeta^*$ can be translated as follows:
\begin{itemize}
\item For every $x\in X$, there exists $y\in Y$ such that $y$ is the image of $x$ under $\zeta$; because $\id_X\leq\zeta^*\circ\zeta$ means that
\begin{equation} \label{idX-leq-zeta}
k\leq\bv\limits_{y\in Y}\zeta^*(y,x)\with\zeta(x,y)
\end{equation}
for all $x\in X$.
\item If $y,z\in Y$ are both the images of $x$ under $\zeta$, then $y$ is equal to $z$; because $\zeta\circ\zeta^*\leq\id_Y$ means that
\begin{equation} \label{zeta-leq-idY}
\bv_{x\in X}\zeta(x,z)\with\zeta^*(y,x)\leq\id_Y(y,z)
\end{equation}
for all $y,z\in Y$.
\end{itemize}
\end{rem}

\begin{rem} \label{Q-map-distinguish}
The notion of $\sQ$-map should be carefully distinguished from related notions in the literature:
\begin{enumerate}[label=(\arabic*)]
\item \emph{Fuzzy mappings} in \cite{Chang1972} and \emph{fuzzifying functions} in \cite{Dubois1980} are precisely our $\sQ$-relations.
\item When $\sQ$ is an integral quantale, every $\sQ$-map is a \emph{fuzzy function} with respect to $\id_X$ and $\id_Y$ in the sense of \cite{Sasaki1993,Demirci1999} or \cite{Demirci2000,Demirci2004}, but not vice versa.
\item \label{Q-map-distinguish:Nemitz} When $\sQ$ is a frame, a $\sQ$-map $\zeta\colon X\rto Y$ with $\zeta\circ\zeta^*=\id_Y$ is exactly a \emph{$\sQ$-valued fuzzy function} in the sense of \cite{Nemitz1986}.
\end{enumerate}
\end{rem}

\begin{lem} \label{zeta-zeta*-k}
Let $\zeta\colon X\rto Y$ be a $\sQ$-map. Then
\[(\zeta^*\circ\zeta)(x,x)=k\quad\text{and}\quad\zeta(x,z)\with\zeta^*(y,x)=\bot\]
for all $x\in X$, $y,z\in Y$ with $y\neq z$.
\end{lem}

\begin{proof}
The first equality follows from
\begin{align*}
k&=\id_X(x,x)\leq(\zeta^*\circ\zeta)(x,x)=\bv_{y\in Y}\zeta^*(y,x)\with\zeta(x,y)\\
&=\bv_{y\in Y}\zeta(x,y)\with\zeta^*(y,x)\leq\bv_{y\in Y}(\zeta\circ\zeta^*)(y,y)\leq\bv_{y\in Y}\id_Y(y,y)=k,
\end{align*}
and the second equality holds because
\[\zeta(x,z)\with\zeta^*(y,x)\leq(\zeta\circ\zeta^*)(y,z)\leq\id_Y(y,z)=\bot.\qedhere\]
\end{proof}

Given a $\sQ$-map $\zeta\colon X\rto Y$, it is obvious that the restriction of $\zeta$ on any (crisp) subset $W\subseteq X$ induces a $\sQ$-map $\zeta|W\colon W\rto Y$. Moreover, $\zeta$ can also be extended to any (crisp) superset $Z\supseteq Y$:

\begin{lem} \label{Qmap-extend}
Let $\zeta\colon X\rto Y$ be a $\sQ$-map, and let $Z\supseteq Y$. Then
\[\tilde{\zeta}\colon X\rto Z,\quad \tilde{\zeta}(x,z)=\begin{cases}
\zeta(x,z) & \text{if}\ z\in Y,\\
\bot & \text{if}\ z\not\in Y
\end{cases}\]
is the unique $\sQ$-map from $X$ to $Z$ satisfying 
\begin{equation} \label{restrict}
\tilde{\zeta}(x,y)=\zeta(x,y)\quad\text{and}\quad(\tilde{\zeta})^*(y,x)=\zeta^*(y,x)
\end{equation}
for all $x\in X$, $y\in Y$.
\end{lem}

\begin{proof}
It is easy to see that $\tilde{\zeta}$ and
\[\eta\colon Z\rto X,\quad\eta(z,x)=\begin{cases}
\zeta^*(z,x) & \text{if}\ z\in Y,\\
\bot & \text{if}\ z\not\in Y
\end{cases}\]
satisfy \eqref{idX-leq-zeta} and \eqref{zeta-leq-idY}, and thus $\tilde{\zeta}\dv\eta$, making $\tilde{\zeta}$ a $\sQ$-map. For the uniqueness of $\tilde{\zeta}$, let $\zeta'\colon X\rto Z$ be another $\sQ$-map satisfying \eqref{restrict}. Then for any $x\in X$ and $z\in Z\setminus Y$,
\begin{align*}
\zeta'(x,z)&\leq\zeta'(x,z)\with\Big(\bv\limits_{y\in Y}\zeta^*(y,x)\with\zeta(x,y)\Big)&(\text{Inequality \eqref{idX-leq-zeta}})\\
&=\bv\limits_{y\in Y}\zeta'(x,z)\with\zeta^*(y,x)\with\zeta(x,y)\\
&=\bv\limits_{y\in Y}\zeta'(x,z)\with\zeta'^*(y,x)\with\zeta(x,y)&(\text{Equations \eqref{restrict}})\\
&\leq\bv\limits_{y\in Y}\id_Z(y,z)\with\zeta(x,y)&(\text{Inequality \eqref{zeta-leq-idY}})\\
&=\bv\limits_{y\in Y}\bot\with\zeta(x,y)\\
&=\bot,
\end{align*}
which shows that $\zeta'=\tilde{\zeta}$.
\end{proof}

Given (crisp) sets $X$, $Y$, we say that a $\sQ$-map $\zeta\colon X\rto Y$ is \emph{symmetric} \cite{Heymans2010} if $\zeta^*=\zeta^{\op}$, i.e., if
\[\zeta^*(y,x)=\zeta(x,y)\]
for all $x\in X$, $y\in Y$. In particular, every map $f\colon X\to Y$ between (crisp) sets induces a symmetric $\sQ$-map
\[f_{\circ}\colon X\rto Y,\quad f_{\circ}(x,y)=\begin{cases}
k & \text{if}\ y=f(x),\\
\bot & \text{else},
\end{cases}\]
called the \emph{graph} of $f$. Thus we have a faithful functor
\[(-)_{\circ}\colon\Set\to\QMap.\]

\begin{lem} \label{zeta_s}
Let $\zeta\colon X\rto Y$ be a $\sQ$-map. Then $\zeta$ is symmetric if, and only if, the $\sQ$-relation
\[\zeta_{\sfs}\colon X\rto Y,\quad\zeta_{\sfs}(x,y)=\zeta(x,y)\wedge\zeta^*(y,x)\]
is a symmetric $\sQ$-map.
\end{lem}

\begin{proof}
If $\zeta$ is symmetric, then $\zeta^*(y,x)=\zeta^{\op}(y,x)=\zeta(x,y)$, and consequently $\zeta_{\sfs}(x,y)=\zeta(x,y)$ for all $x\in X$, $y\in Y$. Thus $\zeta_{\sfs}=\zeta$ is a symmetric $\sQ$-map. Conversely, suppose that $\zeta_{\sfs}$ is a symmetric $\sQ$-map. It follows that $\zeta_{\sfs}\leq\zeta$ and $\zeta_\sfs^\op\leq\zeta^*$ from
\[\zeta_\sfs(x,y)=\zeta(x,y)\wedge\zeta^*(y,x)\leq \zeta(x,y)\]
and 
\[\zeta_\sfs^\op(y,x)=\zeta_\sfs(x,y)=\zeta(x,y)\wedge\zeta^*(y,x)\leq \zeta^*(y,x)\]
for any $x\in X$, $y\in Y$. Then
\begin{align*}
\zeta&=\zeta\circ\id_X\\
&\leq \zeta\circ\zeta_\sfs^*\circ\zeta_\sfs &(\zeta_{\sfs}\dv\zeta_{\sfs}^*)\\
&=\zeta\circ\zeta_\sfs^\op\circ\zeta_\sfs &(\zeta_{\sfs}\ \text{is symmetric})\\
&\leq \zeta\circ\zeta^*\circ\zeta_\sfs\\
&\leq\id_Y\circ\zeta_\sfs&(\zeta\dv\zeta^*)\\
&=\zeta_\sfs\\
&\leq \zeta.
\end{align*}
Thus $\zeta=\zeta_{\sfs}$ is a symmetric $\sQ$-map.    
\end{proof}

The aim of this section is to answer the following questions:
\begin{itemize}
\item Is every $\sQ$-map symmetric?
\item Is every $\sQ$-map the graph of a map in $\Set$?
\end{itemize}
In fact, both the answers are negative, and we will provide necessary and sufficient conditions on $\sQ$ for them to be true.

\begin{defn} \label{lean-def}
Let $\sQ$ be a non-trivial, commutative and unital quantale. We say that:
\begin{enumerate}[label=(\arabic*)]
\item \label{lean-def:lean} $\sQ$ is \emph{lean}, if 
\begin{equation} \label{lean-def:lean:p-vee-q=k}
(p\vee q=k\ \text{and}\ p\with q=\bot)\implies(p=k\ \text{or}\ q=k)
\end{equation}
and
\begin{equation} \label{lean-def:lean:p-with-q=k}
p\with q=k\iff p=q=k
\end{equation}
for all $p,q\in\sQ$;
\item \label{lean-def:weakly-lean} $\sQ$ is \emph{weakly lean}, if
\begin{equation} \label{lean-def:weakly-lean:bv}
\Big(\dbv\limits_{i\in I}p_i\with q_i=k\ \text{and}\ p_i\with q_j=\bot \ (i\neq j)\Big)\implies\Big(k\leq \dbv\limits_{i\in I}(p_i\wedge q_i)\with(p_i\wedge q_i)\Big)
\end{equation}
for all $p_i,q_i\in\sQ$ $(i\in I)$.
\end{enumerate}
\end{defn}

\begin{rem} \label{lean-origin}
\begin{enumerate}[label=(\arabic*)]
\item \label{lean-origin:lean} Note that \eqref{lean-def:lean:p-with-q=k} is always true when $\sQ$ is integral. Indeed, $p\with q=k$ implies that
\[k=p\with q\leq p\with k=p\leq k,\]
which forces $p=k$, and $q=k$ follows from similar calculations. Hence, Definition \ref{lean-def}\ref{lean-def:lean} generalizes the notion of \emph{lean} in \cite[Section III.1.2]{Hofmann2014}, where the quantale $\sQ$ is required to be integral, and consequently, only \eqref{lean-def:lean:p-vee-q=k} is presented there.
\item \label{lean-origin:weakly-lean} The condition \eqref{lean-def:lean:p-with-q=k} is a modification of the condition given in \cite[Proposition 3.5.3(3)]{Heymans2010} and \cite[Theorem 3.7(5)]{Heymans2011}. The difference is that we deal with the discrete case, i.e., $\sQ$-maps between (crisp) sets (instead of \emph{$\sQ$-categories} as in \cite{Heymans2010,Heymans2011}).
\end{enumerate}
\end{rem}

\begin{lem} \label{pvq=k-pwithq=bot-idem}
For $p,q\in\sQ$, if $p\vee q=k$ and $p\with q=\bot$, then both $p$ and $q$ are idempotent.
\end{lem}

\begin{proof}
Just note that
\[p=p\with k=p\with (p\vee q)=(p\with p)\vee(p\with q)=(p\with p)\vee\bot=p\with p,\]
and $q=q\with q$ can be obtained similarly.
\end{proof}


\begin{lem} \label{lean-imp-weakly-lean}
If $\sQ$ is lean, then $\sQ$ is weakly lean.
\end{lem}

\begin{proof}
Suppose that $p_i,q_i\in\sQ$ $(i\in I)$, $\dbv\limits_{i\in I}p_i\with q_i=k$ and $p_i\with q_j=\bot$ $(i\neq j)$. Since $\bot<k$, there exists $j\in I$ such that $p_{j}\with q_{j}>\bot$. Let
\[u:=p_j\with q_j\quad\text{and}\quad v:=\bv_{\substack{i\in I\\ i\neq j}}p_j\with q_j.\]
Then 
\begin{equation} \label{u-v-k-bot}
u\vee v=k\quad\text{and}\quad u\with v=p_{j}\with q_{j}\with\Big(\bv_{\substack{i\in I\\ i\neq j}}p_i\with q_i\Big)=\bv_{\substack{i\in I\\ i\neq j}}p_{j}\with q_{j}\with p_i\with q_i=\bot.
\end{equation}
Note that $v\neq k$. Indeed, if $v=k$, then $u\with v=u\with k=u>\bot$, contradicting to \eqref{u-v-k-bot}. Thus, since $\sQ$ is lean, we deduce from \eqref{u-v-k-bot} that $u=p_{j}\with q_{j}=k$, and consequently $p_{j}=q_{j}=k$. It follows that
\[\bv\limits_{i\in I}(p_i\wedge q_i)\with(p_i\wedge q_i)\geq(p_{j}\wedge q_{j})\with(p_{j}\wedge q_{j})=(k\wedge k)\with(k\wedge k)=k,\]
showing that $\sQ$ is weakly lean.
\end{proof}

\begin{lem} \label{divisible-weakly-lean}
If $\sQ$ is integral, then $\sQ$ is weakly lean. 
\end{lem}

\begin{proof}
Suppose that $p_i,q_i\in\sQ$ $(i\in I)$, $\dbv\limits_{i\in I}p_i\with q_i=k$ and $p_i\with q_j=\bot$ $(i\neq j)$. Note that for each $i\in I$,
\begin{align*}
p_i&=p_i\with k\\
&=p_i\with\Big(\bv_{j\in I}p_j\with q_j\Big)\\
&=\bv_{j\in I}p_i\with p_j\with q_j\\
&=(p_i\with p_i\with q_i)\vee\Big(\bv_{\substack{j\in I\\j\neq i}}p_i\with p_j\with q_j\Big)\\
&=(p_i\with p_i\with q_i)\vee\bot\\
&=p_i\with p_i\with q_i.
\end{align*}
Since $\sQ$ is integral, we have $p\with q\leq (p\with k)\wedge(k\with q)=p\wedge q$ for all $p,q\in\sQ$, and consequently
\[\bv\limits_{i\in I}(p_i\wedge q_i)\with(p_i\wedge q_i)\geq\bv\limits_{i\in I}p_i\with q_i\with p_i\with q_i=\bv\limits_{i\in I}p_i\with q_i=k,\]
as desired.
\end{proof}

\begin{exmp} \label{lean-exmp}
For the examples listed in \ref{quantale-example}:
\begin{enumerate}[label=(\arabic*)]
\item $(C_3,\with,k)$ is a non-integral lean quantale.
\item Every frame is divisible, and thus integral. Then it follows from Lemma \ref{divisible-weakly-lean} that every frame is weakly lean. Moreover: 
\begin{itemize}
\item $F_1=\{\bot,p,q,\top\}$ is not lean, while $F_2=\{\bot,p,q,r,\top\}$ is lean.
\item The powerset $\PX$ is lean if, and only if, $X$ is a singleton set.
\end{itemize}
Therefore, a weakly lean quantale need not be lean.
\item \label{lean-exmp:M3} The quantale $M_3=(M_3,\with,k)$ is not weakly lean, while $M'_3=(M_3,\with',k)$ is lean. 
\item For each continuous t-norm $*$ on $[0,1]$, the quantale $([0,1],*,1)$ is lean. 
\item The quantale $([-\infty,\infty],+,0)$ is not weakly lean, while the Lawvere quantale $([0,\infty],+,0)$ is lean.
\item Consider the free quantale $(\PM,\with,\{k\})$ induced by a commutative monoid $(M,\with,k)$:
\begin{itemize}
\item $(\PM,\with,\{k\})$ is lean if, and only if, $k$ is the only element of $M$ with an inverse; that is, if $m\in M$ and $m\neq k$, then there exists no $m'\in M$ with $m\with m'=k$.
\item $(\PM,\with,\{k\})$ is weakly lean if, and only if, there exist no $m,m'\in M$ such that $m\neq m'$ and $m\with m'=k$.
\end{itemize}
In particular, the free quantale induced by the cyclic group $\bbZ_2=\bbZ/2\bbZ=\{0,1\}$ (under the usual addition modulo $2$) is non-integral, weakly lean, but not lean.
\end{enumerate}
\end{exmp}

The following theorem is a discrete counterpart of \cite[Proposition 3.5.3]{Heymans2010} and \cite[Theorem 3.7]{Heymans2011}, and can be proved in a similar manner:

\begin{thm} \label{Q-map-symmetric-weakly-lean}
Let $\sQ$ be a non-trivial, commutative and unital quantale. Then every $\sQ$-map is symmetric if, and only if, $\sQ$ is weakly lean.
\end{thm}

\begin{proof}
For the necessity, suppose that every $\sQ$-map is symmetric. Assume that $p_i,q_i\in\sQ$ $(i\in I)$, $\dbv\limits_{i\in I}p_i\with q_i=k$ and $p_j\with q_i=\bot$ $(i\neq j)$. Define $\sQ$-relations
\[\zeta\colon\{\star\}\rto I\quad \text{and}\quad \eta\colon I\rto \{\star\}\]
by
\[\zeta(\star,i)=p_i\quad\text{and}\quad \eta(i,\star)=q_i\]
for any $i\in I$ where $\{\star\}$
refers to the singleton set. Note that for each $i\in I$ and $i\neq j$, from $p_i\with q_i\leq k$ and $p_j\with q_i=\bot$ we have
\[(\zeta\circ\eta)(i,i)=\zeta(\star,i)\with\eta(i,\star)=p_i\with q_i\leq k=\id_I(i,i)\]
and
\[(\zeta\circ\eta)(i,j)=\zeta(\star,j)\with\eta(i,\star)=p_j\with q_i=\bot=\id_I(i,j).\]
Thus $\zeta\circ\eta\leq\id_I$. It follows that
\[(\eta\circ\zeta)(\star,\star)=\bv\limits_{i\in I}\eta(i,\star)\with\zeta(\star,i)=\bv\limits_{i\in I}q_i\with p_i=k=\id_{\{\star\}}(\star,\star),\]
which forces that $\zeta\dv\eta$, i.e. $\zeta$ is a $\sQ$-map. Since every $\sQ$-map is symmetric, in particular we have 
\[q_i=\eta(i,\star)=\zeta^\op(i,\star)=\zeta(\star,i)=p_i\]
for all $i\in I$, and consequently
\[\bv\limits_{i\in I}(p_i\wedge q_i)\with (p_i\wedge q_i)=\bv\limits_{i\in I}(p_i\wedge p_i)\with(q_i\wedge q_i)=\bv\limits_{i\in I}p_i\with q_i=k.\]
Hence $\sQ$ is weakly lean.

Conversely, for the sufficiency, suppose that $\sQ$ is 
weakly lean. In order to show that every $\sQ$-map $\zeta\colon X\rto Y$ is symmetric, by Lemma \ref{zeta_s} it suffices to show that $\zeta_{\sfs}$ is always a symmetric $\sQ$-map. To this end, let $x\in X$. 
From Lemma \ref{zeta-zeta*-k} we have
\begin{equation} \label{zeta-zeta*-k-bot}
\bv\limits_{y\in Y}\zeta^*(y,x)\with\zeta(x,y)=k\quad\text{and}\quad\zeta(x,z)\with\zeta^*(y,x)=\bot\ (y\neq z).
\end{equation}
Since $\sQ$ is weakly lean, \eqref{zeta-zeta*-k-bot} forces
\[(\zeta_{\sfs}^\op\circ\zeta_{\sfs})(x,x)=\bv\limits_{y\in Y}\zeta_{\sfs}^{\op}(y,x)\with\zeta_{\sfs}(x,y)=\bv\limits_{y\in Y}(\zeta(x,y)\wedge\zeta^*(y,x))\with(\zeta(x,y)\wedge\zeta^*(y,x))\geq k=\id_X(x,x).\]
Moreover,
\[(\zeta_{\sfs}\circ\zeta_{\sfs}^\op)(y,z)=\bv\limits_{x\in X}\zeta_{\sfs}(x,z)\with\zeta_{\sfs}^{\op}(y,x)\leq \bv\limits_{x\in X}\zeta(x,z)\with\zeta^*(y,x)\leq \id_Y(y,z)\]
for all $y,z\in Y$. Thus $\zeta_{\sfs}\dv\zeta_{\sfs}^{\op}$; that is, $\zeta_{\sfs}$ is a symmetric $\sQ$-map.
\end{proof}

\begin{cor} \label{weakly-lean-map-discrete}
If $\sQ$ is weakly lean, then for any $\sQ$-maps $\zeta,\eta\colon X\rto Y$, $\zeta=\eta$ whenever $\zeta\leq\eta$. 
\end{cor}

\begin{proof}
First, $\zeta\leq\eta$ implies that $\eta^*\leq \zeta^*$. Second, since $\sQ$ is weakly lean, Theorem \ref{Q-map-symmetric-weakly-lean} indicates that $\eta^\op=\eta^*\leq\zeta^*=\zeta^{\op}$. Thus $\eta\leq\zeta$, as desired.
\end{proof}

\begin{rem}
Following the terminology of \cite[Definition 3.1(2)]{Heymans2012a}, Corollary \ref{weakly-lean-map-discrete} actually says that $\QMap$ is \emph{map-discrete} when $\sQ$ is weakly lean.

We would also like to point out that every $\sQ$-map is symmetric if $\sQ$ is \emph{modular} (cf. \cite[Definition 2.5.1]{Heymans2010}) and \emph{localic} (cf. \cite[Definition 2.1.1]{Heymans2010}), which is an immediate consequence of \cite[Propositions 2.5.3(4) and 2.5.9]{Heymans2010}. However, this is only a sufficient condition, since a weakly lean quantale need not be localic; for example, the quantale $M'_3=(M_3,\with',k)$ given in Examples \ref{quantale-example}\ref{quantale-example:diamond} and \ref{lean-exmp}\ref{lean-exmp:M3} is lean but not localic. 
\end{rem}

\begin{thm} \label{QMap-Set-lean}
Let $\sQ$ be a non-trivial, commutative and unital quantale. Then every $\sQ$-map is the graph of a map in $\Set$ if, and only if, $\sQ$ is lean. In this case, $\Set$ and $\QMap$ are isomorphic categories.
\end{thm}

\begin{proof}
For the necessity, suppose that every $\sQ$-map is the graph of a map in $\Set$. 

First, let $p,q\in\sQ$ with $p\vee q=k$ and $p\with q=\bot$. Then it follows from Lemma \ref{pvq=k-pwithq=bot-idem} that both $p$ and $q$ are idempotent. Define $\sQ$-relations
\[\zeta\colon\{\star\}\rto\{p,q\}\quad\text{and}\quad \eta\colon\{p,q\}\rto\{\star\}\]
by
\[\zeta(\star,p)=\eta(p,\star)=p\quad\text{and}\quad \zeta(\star,q)=\eta(q,\star)=q.\]
Then, from
\[\eta\circ\zeta=\begin{pmatrix}
p & q
\end{pmatrix}
\begin{pmatrix}
p\\
q
\end{pmatrix}
=(p\with p)\vee(q\with q)=p\vee q=k=\id_{\{\star\}}\]
and
\[\zeta\circ\eta=\begin{pmatrix}
p\\
q
\end{pmatrix}
\begin{pmatrix}
p & q
\end{pmatrix}=\begin{pmatrix}
p\with p & p\with q\\
q\with p & q\with q
\end{pmatrix}=\begin{pmatrix}
p & \bot\\
\bot & q
\end{pmatrix}\leq\begin{pmatrix}
k & \bot\\
\bot & k
\end{pmatrix}=\id_{\{p,q\}}\]
we obtain that $\zeta\dv\eta$. Thus $\zeta$ is a $\sQ$-map, and consequently it is the graph of a map $f\colon\{\star\}\to\{p,q\}$ in $\Set$, which forces $p=k$ or $q=k$.

Second, let $p,q\in\sQ$ with $p\with q=k$. Define $\sQ$-relations 
\[\zeta_p\colon\{\star\}\rto\{\star\}\quad\text{and}\quad \zeta_q\colon\{\star\}\oto\{\star\}\]
by
\[\zeta_p(\star,\star)=p\quad\text{and}\quad\zeta_q(\star,\star)=q.\]
Then, from $p\with q=k$ we see that $\zeta_p\dv\zeta_q\dv\zeta_p$ trivially holds. Thus, both $\zeta_p$ and $\eta_q$ are $\sQ$-maps, which are necessarily the graph of the unique (identity) map $g\colon\{\star\}\to\{\star\}$. Hence $p=q=k$.

Conversely, for the sufficiency, suppose that $\sQ$ is lean. Let $\zeta\colon X\rto Y$ be a $\sQ$-map. For each $x\in X$, by Lemma \ref{zeta-zeta*-k} we have
$\dbv\limits_{y\in Y}\zeta^*(y,x)\with\zeta(x,y)=k$. Since $k>\bot$, there exists $f(x)\in Y$ with $p:=\zeta^*(f(x),x)\with \zeta(x,f(x))>\bot$. Let
\[q:=\bv\limits_{\substack{y\in Y\\ y\neq f(x)}}\zeta^*(y,x)\with\zeta(x,y).\]
Then $p\vee q=k$, and
\begin{align*}
p\with q&=\zeta^*(f(x),x)\with\zeta(x,f(x))\with \Big(\bv\limits_{\substack{y\in Y\\ y\neq f(x)}}\zeta^*(y,x)\with \zeta(x,y)\Big)\\
&=\bv\limits_{\substack{y\in Y\\ y\neq f(x)}}\zeta^*(f(x),x)\with \zeta(x,f(x))\with \zeta^*(y,x)\with \zeta(x,y)\\
&=\bv\limits_{\substack{y\in Y\\ y\neq f(x)}}\zeta(x,f(x))\with \bot\with \zeta^*(f(x),x))\\
&=\bot,
\end{align*}
where the third equality follows from Lemma \ref{zeta-zeta*-k}.

Note that $q\neq k$, because $q=k$ would lead to $p\with q=p\with k=p>\bot$, contradicting to $p\with q=\bot$. Therefore, given that $\sQ$ is lean, we conclude that $p=\zeta^*(f(x),x)\with \zeta(x,f(x))=k$, and consequently 
\[\zeta^*(f(x),x)=\zeta(x,f(x))=k.\]
We claim that $\zeta$ is the graph of the map $f\colon X\to Y$. Indeed, it is obvious that $f_{\circ}\leq\zeta$, and hence the conclusion follows immediately from Lemma \ref{lean-imp-weakly-lean} and Corollary \ref{weakly-lean-map-discrete}.
\end{proof}

As mentioned in Remark \ref{lean-origin}\ref{lean-origin:lean}, Theorem \ref{QMap-Set-lean} reduces to \cite[Proposition III.1.2.1]{Hofmann2014} when $\sQ$ is integral:

\begin{cor} (See \cite{Hofmann2014}.)
Let $\sQ$ be a non-trivial, commutative and integral quantale. Then every $\sQ$-map is the graph of a map in $\Set$ if, and only if,
\[(p\vee q=\top\ \text{and}\ p\with q=\bot)\implies(p=\top\ \text{or}\ q=\top)\]
for $p,q\in\sQ$.
\end{cor}

As an application of $\sQ$-maps, we extend the notion of frame-valued partition proposed in \cite[Definition 2.11]{Nemitz1986} to the quantale-valued setting:


\begin{defn} \label{Q-partition}
A \emph{$\sQ$-partition} of a (crisp) set $X$ is a collection $\Si$ of $\sQ$-subsets of $X$ (i.e., functions $X\to\sQ$), such that
\begin{enumerate}[label=(P\arabic*)]
\item \label{Q-partition:ST} $Sx\with Tx=\bot$ for all $x\in X$ and $S,T\in\Si$ with $S\neq T$,
\item \label{Q-partition:S} $\dbv\limits_{S\in\Si}Sx=k$ for all $x\in X$,
\item \label{Q-partition:x} $\dbv\limits_{x\in X}Sx=k$ for all $S\in\Si$.
\end{enumerate}
\end{defn}


A $\sQ$-map $\zeta\colon X\rto Y$ is said to be \emph{surjective} if $\zeta\circ\zeta^*=\id_Y$ (cf. Remark \ref{Q-map-distinguish}\ref{Q-map-distinguish:Nemitz}). As the generalization of \cite[Theorem 2.14, Lemma 2.16 and Theorem 2.17]{Nemitz1986}, in what follows we show that every surjective $\sQ$-map gives rise to a $\sQ$-partition, and vice versa. 

For each $\sQ$-map $\zeta\colon X\rto Y$, define
\[\Si_{\zeta}=\{S_y\mid y\in Y\},\]
where
\[S_y\colon X\to\sQ,\quad S_y x=\zeta(x,y)\with \zeta^*(y,x).\]

\begin{prop} \label{Si-zeta}
If $\zeta\colon X\rto Y$ is a surjective $\sQ$-map, then $\Si_{\zeta}$ is a $\sQ$-partition of $X$, and $S_y\neq S_{y'}$ whenever $y\neq y'$ in $Y$.
\end{prop}

\begin{proof}
Suppose that $y\neq y'$ and $S_y=S_{y'}$. For each $x\in X$, note that
\begin{equation} \label{Szx=k}
\bv\limits_{z\in Y}S_z x=\bv\limits_{z\in Y}\zeta^*(z,x)\with\zeta(x,z)=k\quad\text{and}\quad\zeta(x,z)\with\zeta^*(y,x)=\bot\quad (z\neq y)
\end{equation}
by Lemma \ref{zeta-zeta*-k}. Thus, it is easy to deduce that $S_y x=\zeta(x,y)\with\zeta^*(y,x)$ is idempotent by Lemma \ref{pvq=k-pwithq=bot-idem}, and
\begin{equation} \label{Syy'-bot}
S_y x\with S_{y'}x=\zeta(x,y)\with\zeta^*(y,x)\with\zeta(x,y')\with\zeta^*(y',x)=\bot.
\end{equation}
It follows that
\begin{align}
k&=(\zeta\circ\zeta^*)(y,y) & (\zeta\circ\zeta^*=\id_Y) \nonumber\\
&=\bv\limits_{x\in X}S_y x & (S_y x=\zeta(x,y)\with\zeta^*(y,x)) \label{Syx=k}\\
&=\bv\limits_{x\in X}S_y x\with S_y x & (S_y x\ \text{is idempotent}) \nonumber\\
&=\bv\limits_{x\in X}S_y x\with S_{y'}x &(S_y=S_{y'}) \nonumber\\
&=\bot, \nonumber
\end{align}
which is a contradiction. To see that $\Si_{\zeta}$ is a $\sQ$-partition of $X$, just note that \ref{Q-partition:ST}, \ref{Q-partition:S} and \ref{Q-partition:x} are already obtained in \eqref{Syy'-bot}, \eqref{Szx=k} and \eqref{Syx=k}, respectively.
\end{proof}

\begin{prop} \label{zeta-Si}
Let $\Si$ be a $\sQ$-partition of a (crisp) set $X$. Then
\[\zeta_{\Si}\colon X\rto\Si,\quad \zeta_{\Si}(x,S)=Sx\]
is a surjective $\sQ$-map.
\end{prop}

\begin{proof}
From \ref{Q-partition:ST}, \ref{Q-partition:S} and Lemma \ref{pvq=k-pwithq=bot-idem} we obtain that $Sx$ is idempotent for all $x\in X$, $S\in\Si$. It follows that
\[\bv\limits_{S\in\Si}\zeta_{\Si}^\op(S,x)\with\zeta_{\Si}(x,S)=\bv\limits_{S\in\Si}\zeta_{\Si}(x,S)\with\zeta_{\Si}(x,S)=\bv\limits_{S\in\Si}Sx\with Sx=\bv\limits_{S\in\Si}Sx=k,\]
and similarly $\dbv\limits_{x\in X}\zeta_{\Si}(x,S)\with\zeta_{\Si}^\op(S,x)=k$. Moreover,
\[\zeta_{\Si}(x,S)\with \zeta_{\Si}^\op(T,x)=\zeta_{\Si}(x,S)\with\zeta_{\Si}(x,T)=Sx\with Tx=\bot\]
for all $x\in X$ and $S,T\in\Si$ with $S\neq T$. Thus $\zeta_{\Si}\dv\zeta_{\Si}^\op$ and $\zeta_{\Si}\circ\zeta_{\Si}^\op=\id_{\Si}$.
\end{proof}

\begin{prop} \label{Si-zeta-Si=Si}
For each $\sQ$-partition $\Si$ of a (crisp) set $X$, it holds that $\Si_{\zeta_{\Si}}=\Si$.
\end{prop}

\begin{proof}
From the proof of Proposition \ref{zeta-Si} we see that $\zeta_{\Si}\dv\zeta_{\Si}^\op$, and $Tx$ is idempotent for all $x\in X$, $T\in\Si$. Thus $\Si_{\zeta_{\Si}}$ consists of $\sQ$-subsets $S_T\colon X\to\sQ$ of $X$, with
\[S_T x=\zeta_{\Si}(x,T)\with\zeta_{\Si}^{\op}(T,x)=Tx\with Tx=Tx\]
for all $T\in\Si$, $x\in X$. Hence $\Si_{\zeta_{\Si}}=\Si$.
\end{proof}

\begin{exmp}
We have seen several quantales that are not lean in Example \ref{lean-exmp}. In these cases, Theorem \ref{QMap-Set-lean} guarantees the existence of $\sQ$-maps which are not induced by (i.e., the graphs of) maps in $\Set$:
\begin{enumerate}[label=(\arabic*)]
\item Let $X=\{x,y\}$. Then
\[\zeta\colon\{\star\}\rto X,\quad\zeta(\star,x)=p\quad\text{and}\quad \zeta(\star,y)=q\]
is an $F_1$-map, and
\[\eta\colon X\rto X,\quad\eta(x,x)=\eta(y,y)=p\quad \text{and}\quad \eta(x,y)=\eta(y,x)=q.\]
is a surjective $F_1$-map. The $F_1$-partition of $X$ induced by $\eta$ is given by
\[\Si_\eta=\{S_x,S_y\},\quad S_x x=S_y y=p,\quad S_x y=S_y x=q.\]
\item Let $X=\{x,y\}$. Then
\[\zeta\colon \{\star\}\rto X,\quad\zeta(\star,x)=a\quad \text{and}\quad\zeta(\star,y)=\bot\]
is an $M_3$-map, with 
\[\zeta^*\colon X\rto\{\star\},\quad \zeta^*(x,\star)=b\quad \text{and}\quad\zeta^*(y,\star)=\bot.\] 
Note that $\zeta^*\neq\zeta^{\op}$, which may occur because $M_3$ is not weakly lean.
\item Let $X=\{a,b,c\}$, $Y=\{x,y,z\}$ and $Z=\{l,m\}$. Then $\zeta,\eta\colon \{\star\}\rto Y$ given by
\[\begin{array}{c|c}
	       \zeta&\star\\
	\hline x&\{a\}\\
    y&\{b,c\}\\
    z&\varnothing\\
\end{array}\quad\text{and}\quad \begin{array}{c|c}
	       \eta&\star\\
	\hline x&\{a\}\\
    y&\{b\}\\
    z&\{c\}\\
\end{array}\]
are $\PX$-maps, and $\xi\colon Y\rto Z$ given by
\[\begin{array}{c|ccc}
	       \xi&x&y&z\\
	\hline l&\{a\}&\{b\}&\{c\}\\
    m&\{b,c\}&\{a,c\}&\{a,b\}\\
\end{array}\]
is a surjective $\PX$-map. Indeed, the latter follows from
\[\xi\circ\xi^\op=\begin{pmatrix}
\{a\} & \{b\} &\{c\}\\
\{b,c\}& \{a,c\}& \{a,b\}
\end{pmatrix}\begin{pmatrix}
\{a\}& \{b,c\}\\
\{b\}&\{a,c\}\\
\{c\}&\{a,b\}
\end{pmatrix}=\begin{pmatrix}
X & \varnothing\\
\varnothing& X
\end{pmatrix}=\id_Z\]
and
\[\xi^\op\circ\xi=\begin{pmatrix}
\{a\}& \{b,c\}\\
\{b\}&\{a,c\}\\
\{c\}&\{a,b\}
\end{pmatrix}\begin{pmatrix}
\{a\} & \{b\} &\{c\}\\
\{b,c\}& \{a,c\}& \{a,b\}
\end{pmatrix}=\begin{pmatrix}
X & \{c\} &\{b\}\\
\{c\}& X& \{a\}\\
\{b\}&\{a\}& X
\end{pmatrix}\geq \id_Y.\]
The $\PX$-partition $\Si_{\xi}=\{S_l,S_m\}$ of $Y$ induced by $\xi$ is given by
\[\begin{array}{c|cc}
   &S_l  &S_m  \\
\hline   x  & \{a\} & \{b,c\}\\
y & \{b\} & \{a,c\}\\
z&\{c\} &\{a,b\}
\end{array}.\]
\item Let $X=\{x,y\}$. Then $\zeta,\eta\colon X\rto X$ given by
\[\begin{array}{c|cc}
  \zeta &x &y  \\
\hline   x  & -1 & \infty\\
y & \infty & 1
\end{array}\quad \text{and}\quad \begin{array}{c|cc}
  \eta &x &y  \\
\hline   x  & 1 & \infty\\
y & \infty & -1
\end{array}\]
are surjective $[-\infty,\infty]$-maps, because
\[\eta\circ\zeta=\begin{pmatrix}
1&\infty\\
\infty& -1
\end{pmatrix}\begin{pmatrix}
-1&\infty\\
\infty& 1
\end{pmatrix}=\begin{pmatrix}
0&\infty\\
\infty& 0
\end{pmatrix}=\id_X\]
and
\[\zeta\circ\eta=\begin{pmatrix}
-1&\infty\\
\infty& 1
\end{pmatrix}\begin{pmatrix}
1&\infty\\
\infty& -1
\end{pmatrix}=\begin{pmatrix}
0&\infty\\
\infty& 0
\end{pmatrix}=\id_X.\]
The $[-\infty,\infty]$-partition $\Si_{\zeta}=\{S_x,S_y\}$ of $X$ induced by $\zeta$ is given by
\[\begin{array}{c|cc}
   &S_x  &S_y  \\
\hline   x  & 0 & \infty\\
y & \infty & 0
\end{array},\]
which can be regarded as the crisp partition $\{\{x\},\{y\}\}$ of $X$.
\end{enumerate}
\end{exmp}

\section{Quantale-valued partial maps} \label{Quantale-valued_partial_maps}

Note that coproducts of any family $\{X_i\}_{i\in I}$ of sets in the category $\QMap$ are given by disjoint unions $X=\coprod\limits_{i\in I}X_i$ of sets as in $\Set$, and the graphs $\iota_i\colon X_i\rto X$ $(i\in I)$ of the inclusion maps are coproduct injections.



\begin{lem} \label{zeta-amalg-eta}
If $\zeta\colon X\rto Y$, $\eta\colon X'\rto Y'$ are $\sQ$-maps, then 
\[\zeta\amalg\eta\colon X\amalg X'\rto Y\amalg Y'\] 
is a $\sQ$-map. 
\end{lem}

\begin{proof}
It is straightforward to check that $\zeta^*\amalg\eta^*\colon Y\amalg Y'\rto X\amalg X'$ is the right adjoint of $\zeta\amalg\eta$ in $\QRel$.
\end{proof}

Let
\[X_+:=X\amalg\{\star\}\]
denote the disjoint union of every (crisp) set $X$ and the singleton set $\{\star\}$.

\begin{defn}
Let $X$, $Y$ be (crisp) sets. A \emph{partial $\sQ$-map} $\zeta$ from $X$ to $Y$ is a $\sQ$-map $\zeta\colon X\rto Y_+$.
\end{defn}

\begin{rem}
The singleton set $\{\star\}$ may not be a terminal object of $\QMap$. For example, if $\sQ=M_3$ (see Example \ref{quantale-example}\ref{quantale-example:diamond}), then there are two different $\sQ$-maps from $\{\star\}$ to $\{\star\}$:
\begin{itemize}
\item the graph of the identity map on $\{\star\}$;
\item the $\sQ$-map
\[\zeta\colon\{\star\}\rto\{\star\},\quad \zeta(\star,\star)=a,\]
whose right adjoint is given by
\[\eta\colon\{\star\}\rto\{\star\},\quad \eta(\star,\star)=b.\]
\end{itemize}
\end{rem}

\begin{rem}
When $\sQ$ is an integral quantale, every partial $\sQ$-map is a \emph{fuzzy partial function} defined in \cite{Demirci2000,Demirci2004}, but not vice versa.
\end{rem}

In order to construct the category of sets and partial $\sQ$-maps, note that by setting 
\[\CC=\QMap\quad\text{and}\quad A=\{\star\}\]
in \eqref{U-AC-C} and \eqref{A-AC} we obtain a pair of adjoint functors
\begin{equation} \label{star-QMap-adjoint}
\bfig
\morphism|a|/@{->}@<6pt>/<1000,0>[\QMap`\sQMap;\sF]
\morphism(1000,0)|b|/@{->}@<6pt>/<-1000,0>[\sQMap`\QMap;\sU]
\place(450,0)[\bot]
\efig
\end{equation}
given by the following data:
\begin{itemize}
\item $\sU$ is the forgetful functor;
\item the functor $\sF$ sends each set $X$ to the graph of the inclusion map $\{\star\}\ \to/^(->/X_+$, and sends each $\sQ$-map $\zeta\colon X\rto Y$ to the $\sQ$-map 
\[\zeta_+:=\zeta\amalg\id_{\{\star\}}\colon X_+\rto Y_+;\]
\item for each set $X$, the component of the unit $\iota_X \colon X\rto X_+$ is the graph of the inclusion map $X\ \to/^(->/X_+$; 
\item for each object $\mu\colon \{\star\}\rto X$ in $\sQMap$, the component of the counit $\ep_{\mu} \colon \sF X \rto \mu$ is given by 
\[\ep_{\mu}(x,y)=\begin{cases}
\id_X(x,y) & \text{if}\ x,y\in X,\\
\mu(\star,y) & \text{if}\ x=\star\ \text{and}\ y\in X.
\end{cases}\]
\end{itemize}
The adjunction $\sF\dashv\sU$ induces the \emph{maybe monad} $(\sT,\sm,\iota)$ on $\QMap$, whose endofunctor $\sT=\sU\sF$ carries each set $X$ to $X_+$, and whose multiplication is given by
\begin{equation} \label{m_X-def}
\sm_{X}\colon (X_+)_+\rto X_+,\quad\sm_X(x,y)=\begin{cases}
\id_X(x,y) & \text{if}\ x,y\in X,\\
k & \text{if}\ x=y=\star,\\
\bot & \text{if}\ x=\star,\ y\in X\ \text{or}\ x\in X,\ y=\star.
\end{cases}
\end{equation}
The isomorphism between $\sQMap$ and the Eilenberg-Moore category $\QMap^{\sT}$ is not difficult to be observed. Indeed, a $\sT$-algebra $(X,\mu)$ consists of a set $X$ and a $\sQ$-map $\mu\colon X_+\rto X$ such that 
\[\mu\circ\iota_X=\id_X,\]
which is clearly uniquely determined by an object in $\sQMap$. Moreover, it is easy to see that a morphism $\zeta\colon (X,\mu)\rto(Y,\lambda)$ of $\sT$-algebras is a $\sQ$-map $\zeta\colon X\rto Y$ such that 
\[(\zeta\circ\mu)(\star,y)=\lambda(\star,y)\]
for all $y\in Y$, which is simply a morphism in $\sQMap$. Therefore, in what follows we identify $\QMap^{\sT}$ with $\sQMap$.


Furthermore, objects in the Kleisli category $\QMap_{\sT}$ are sets, and a morphism from $X$ to $Y$ in $\QMap_{\sT}$ is exactly a $\sQ$-map $X\rto Y_+$; that is, a partial $\sQ$-map from $X$ to $Y$. The composite of partial $\sQ$-maps $\zeta\colon X\rto Y_+$ and $\eta\colon Y\rto Z_+$ is given by 
\begin{equation} \label{QParMap-comp}
(\eta\diamond\zeta)(x,z)=\begin{cases}
\dbv\limits_{y\in Y} \eta(y,z) \with \zeta(x,y) & \text{if}\ x\in X\ \text{and}\ z\in Z,\\
\Big(\dbv\limits_{y\in Y}\eta(y,\star)\with \zeta(x,y)\Big)\vee\zeta(x,\star) & \text{if}\ x\in X\ \text{and}\ z=\star.
\end{cases}
\end{equation}
Thus, we denote by
\[\QParMap:=\QMap_{\sT}\]
the category of sets and partial $\sQ$-maps. By Lemma \ref{Qmap-extend}, $\QMap$ is embedded in $\QParMap$ as a full subcategory.

\begin{rem} \label{why-partial-Q-map}
Let $\zeta\colon X\to Y_+$ be a partial $\sQ$-map from $X$ to $Y$. The value
\[\zeta(x,\star)\]
is interpreted as the degree that $x$ has no image under the map $\zeta$, and so is the value
\[\zeta^*(\star,x)=\bw_{z\in Y_+}\zeta(x,z)\ra\id_{Y_+}(\star,z),\]
because the above expression means that
\begin{itemize}
\item For each $z\in Y_+$, if $z$ is the image of $x$ under $\zeta$, then $z$ is equal to $\star$; in other words, $x$ has no image under $\zeta$ within $Y$.
\end{itemize}  
The adjunction $\zeta\dashv\zeta^*$ can be translated as follows:
\begin{itemize}
\item For every $x\in X$, there may not be $y\in Y$ such that $y$ is the image of $x$ under $\zeta$; because $\id_X\leq\zeta^*\circ\zeta$ means that
\[k\leq\bv\limits_{y\in Y_+}\zeta^*(y,x)\with\zeta(x,y)=(\zeta^*(\star,x)\with\zeta(x,\star))\vee\Big(\bv\limits_{y\in Y}\zeta^*(y,x)\with\zeta(x,y)\Big)\]
for all $x\in X$.
\item If $y,z\in Y$ are both the images of $x$ under $\zeta$, then $y$ is equal to $z$; because $\zeta\circ\zeta^*\leq\id_{Y_+}$ means that 
\[\bv_{x\in X}\zeta(x,z)\with\zeta^*(y,x)\leq\id_{Y_+}(y,z)=\id_Y(y,z)\]
for all $y,z\in Y$. Moreover, setting $y\in Y$ and $z=\star$ in the above inequality we obtain that
\[\zeta(x,\star)\with\zeta^*(y,x)=\id_{Y_+}(y,\star)=\bot\]
for all $x\in X$; that is, ``$x$ has no image under $\zeta$'' and ``$y$ is the image of $x$ under $\zeta$'' cannot happen simultaneously. 
\end{itemize}
Given another partial $\sQ$-map $\eta\colon Y\rto Z_+$ from $Y$ to $Z$, the formula \eqref{QParMap-comp} is understood as follows:
\begin{itemize}
\item $z\in Z$ is the image of $x$ under $\eta\diamond\zeta$ if, and only if, there exists $y\in Y$ such that $y$ is the image of $x$ under $\zeta$, and $z$ is the image of $y$ under $\eta$.
\item $x$ has no image under $\eta\diamond\zeta$ if, and only if, either of the following cases holds:
    \begin{itemize}
    \item $x$ has no image under $\zeta$;
    \item $x$ has an image $y$ under $\zeta$, but $y$ has no image under $\eta$.
    \end{itemize}
\end{itemize}
\end{rem}

\begin{thm}
Let $\sQ$ be a non-trivial, commutative and unital quantale. Then every partial $\sQ$-map is the graph of a partial map in $\Set$ if, and only if, $\sQ$ is lean. In this case, $\Set^\partial$ and $\QParMap$ are isomorphic categories.
\end{thm}

\begin{proof}
Since $\QMap$ is embedded in $\QParMap$ as a full subcategory, the necessity is an immediate consequence of Theorem \ref{QMap-Set-lean}. For the sufficiency, suppose that $\sQ$ is lean. Since every partial $\sQ$-map is a $\sQ$-map $\zeta\colon X\rto Y_+$, there is a map $f\colon X\to Y_+$ in $\Set$ such that $\zeta=f_{\nat}$, which is exactly a partial map in $\Set$.
\end{proof}

As it is well known that $\Set^{\partial}$ is equivalent to $\{\star\}/\Set$, and thus monadic over $\Set$ \cite{MacLane1998,Riehl2016}, it is natural to ask whether $\QParMap$ has similar conclusions. Explicitly, our question is:
\begin{itemize}
\item Is $\QParMap$ equivalent to $\{\star\}/\QMap$, and thus monadic over $\QMap$?
\end{itemize}
By Beck's monadicity theorem \cite{Beck1967,MacLane1998,Riehl2016}, the answer is affirmative if, and only if, the functor
\[\sU_{\sT}\colon\QParMap\to\QMap,\quad \sU_{\sT} X=X_+,\quad\sU_{\sT}(\zeta\colon X\rto Y_+)=(\sm_Y\circ \sF\zeta\colon X_+\rto Y_+)\]
creates coequalizers of $\sU_{\sT}$-split pairs (see, e.g., \cite[Theorem 5.5.1]{Riehl2016}). However, due to the difficulty of handling coequalizers in $\QMap$ and $\QParMap$, in what follows we proceed in an alternative way.

Let $(X_+,\tau_X)\in\sQMap$ denote the \emph{free} $\sT$-algebra on a set $X$. Then
\begin{equation} \label{tau_X-def}
\tau_X\colon\{\star\}\rto X_+,\quad\tau_X(\star,x)=\begin{cases}
k & \text{if}\ x=\star,\\
\bot & \text{if}\ x\in X
\end{cases}
\end{equation}
is the restriction of $\sm_{X}\colon (X_+)_+\rto X_+$ (see \eqref{m_X-def}), which is exactly the graph of the inclusion map $\{\star\}\ \to/^(->/X_+$. It is well known that the canonical functor
\[\sK\colon\QParMap(=\QMap_{\sT})\to\sQMap(\cong\QMap^{\sT})\]
that sends each $\sQ$-map $\zeta\colon X\rto Y_+$ (i.e., partial $\sQ$-map from $X$ to $Y$) to
\[\sK\zeta\colon (X_+,\tau_X)\rto (Y_+,\tau_Y),\quad (\sK\zeta)(x,y)=\begin{cases}
\zeta(x,y) &\text{if}\ x\in X,\ y\in Y,\\
k &\text{if}\ x=\star=y,\\
\bot &\text{else}
\end{cases}\]
is fully faithful (see, e.g., \cite[Lemma 5.2.13]{Riehl2016}).

\begin{lem} \label{K-sur}
Assuming the axiom of choice, every object of $\sQMap$ is isomorphic to a free $\sT$-algebra, hence the functor $\sK$ is essentially surjective.
\end{lem}

\begin{proof}
Let $(X,\mu)\in\sQMap$. Since $X$ is the codomain of the $\sQ$-map $\mu\colon\{\star\}\rto X$, it is clear that $X\neq\varnothing$. By the axiom of choice and \cite{Hajnal1972}, there exists a binary operation $\cdot$ on $X$ such that $(X,\cdot)$ is a group, with $e\in X$ being its unit. Then
\begin{enumerate}[label=(G\arabic*)]
\item \label{group-mult} $x\cdot z\neq y\cdot z$ and $z\cdot x\neq z\cdot y$ for all $x,y,z\in X$ with $x\neq y$, and
\item \label{group-e} $x\cdot y\neq x$ for all $x,y\in X$ with $y\neq e$.
\end{enumerate}
Our strategy is to show that $(X,\mu)$ is isomorphic to
\[\sK(X\setminus\{e\})=((X\setminus\{e\})_+,\tau_{X\setminus\{e\}}).\]
Explicitly, we define $\sQ$-relations 
\[\zeta\colon X\rto (X\setminus\{e\})_+,\quad\zeta(x,a)=\begin{cases}
\mu^*(x,\star) &\text{if}\ x\in X,\ a=\star\\
\mu^*(x\cdot a,\star) & \text{if}\ x,a\in X
\end{cases}\]
and 
\[\eta\colon (X\setminus\{e\})_+\rto X,\quad \eta(a,x)=\begin{cases}
\mu(\star,x) &\text{if}\ a=\star,\ x\in X\\
\mu(\star,x\cdot a) &\text{if}\ a,x\in X,
\end{cases}\] 
and we show that $\zeta$ and $\eta$ establish an isomorphism between $(X,\mu)$ and $((X\setminus\{e\})_+,\tau_{X\setminus\{e\}})$ in $\sQMap$.

First, $\eta\circ\zeta=\id_X$ and $\zeta\circ\eta=\id_{(X\setminus \{e\})_+}$, so that $\zeta$ and $\eta$ are both $\sQ$-maps, and they are inverses of each other. Indeed, it follows from \ref{group-mult} and Lemma \ref{zeta-zeta*-k} that
\[(\eta\circ\zeta)(x,y)=\bv_{a\in(X\setminus\{e\})_+}\eta(a,y)\with\zeta(x,a)=(\mu(\star,y)\with\mu^*(x,\star))\vee\Big(\bv_{a\in X\setminus\{e\}}\mu(\star,y\cdot a)\with\mu^*(x\cdot a,\star)\Big)=\bot\]
for all $x,y\in X$ with $x\neq y$, and
\begin{align*}
(\eta\circ\zeta)(x,x)&=(\mu(\star,x)\with\mu^*(x,\star))\vee\Big(\bv_{a\in X\setminus\{e\}}\mu(\star,x\cdot a)\with\mu^*(x\cdot a,\star)\Big)\\
&\leq\id_X(x,x)\vee\Big(\bv_{a\in X\setminus\{e\}}\id_X(x\cdot a,x\cdot a)\Big)&(\mu\dv\mu^*)\\
&=k\\
&=(\mu^*\circ\mu)(\star,\star) & (\text{Lemma \ref{zeta-zeta*-k}})\\
&=(\mu(\star,x)\with\mu^*(x,\star))\vee\Big(\bv\limits_{\substack{y\in X\\ y\neq x}}\mu(\star,y)\with\mu^*(y,\star)\Big)\\
&=(\mu(\star,x)\with\mu^*(x,\star))\vee\Big(\bv\limits_{\substack{y\in X\\ y\neq x}}\mu(\star,x\cdot x^{-1}\cdot y)\with\mu^*(x\cdot x^{-1}\cdot y,\star)\Big)\\
&\leq(\mu(\star,x)\with\mu^*(x,\star))\vee\Big(\bv\limits_{a\in X\setminus \{e\}}\mu(\star,x\cdot a)\with\mu^*(x\cdot a,\star)\Big)&(x\neq y\implies x^{-1}\cdot y\neq e)\\
&=(\eta\circ\zeta)(x,x)
\end{align*}
for all $x\in X$. Thus $\eta\circ\zeta=\id_X$. In order to show that $\zeta\circ\eta=\id_{(X\setminus \{e\})_+}$, from \ref{group-mult}, \ref{group-e} and Lemma \ref{zeta-zeta*-k} we deduce that
\begin{align*}
(\zeta\circ\eta)(a,b)&=\bv_{x\in X}\zeta(x,b)\with\eta(a,x)=\bv_{x\in X}\mu^*(x\cdot b,\star)\with\mu(\star,x\cdot a)=\bot,\\
(\zeta\circ\eta)(a,\star)&=\bv_{x\in X}\zeta(x,\star)\with\eta(a,x)=\bv_{x\in X}\mu^*(x,\star)\with\mu(\star,x\cdot a)=\bot,\\
(\zeta\circ\eta)(\star,b)&=\bv_{x\in X}\zeta(x,b)\with\eta(\star,x)=\bv_{x\in X}\mu^*(x\cdot b,\star)\with\mu(\star,x)=\bot
\end{align*}
for all $a,b\in X\setminus \{e\}$ with $a\neq b$. Moreover, it follows from Lemma \ref{zeta-zeta*-k} that
\[(\zeta\circ\eta)(\star,\star)=\bv\limits_{x\in X}\zeta(x,\star)\with\eta(\star,x)=\bv\limits_{x\in X}\mu^*(x,\star)\with\mu(\star,x)=(\mu^*\circ\mu)(\star,\star)=k,\]
and
\begin{align*}
(\zeta\circ\eta)(a,a)&=\bv\limits_{x\in X}\mu(\star,x\cdot a)\with\mu^*(x\cdot a,\star)\\
&\leq\bv\limits_{x\in X}\id_X(x\cdot a,x\cdot a)&(\mu\dv\mu^*)\\
&=k \\
&=(\mu^*\circ\mu)(\star,\star)&(\text{Lemma \ref{zeta-zeta*-k}})\\
&=\bv\limits_{x\in X}\mu^*(x,\star)\with\mu(\star,x)\\
&=\bv\limits_{x\in X}\mu^*(x\cdot a^{-1}\cdot a,\star)\with\mu(\star,x\cdot a^{-1}\cdot a)\\
&\leq\bv\limits_{y\in X}\mu^*(y\cdot a,\star)\with\mu(\star,y\cdot a)\\
&=(\zeta\circ\eta)(a,a)
\end{align*}
for all $a\in X\setminus \{e\}$, as desired.

Second, $\zeta\circ\mu=\tau_{X\setminus \{e\}}$ and $\eta\circ\tau_{X\setminus\{e\}}=\mu$, so that $\zeta\colon(X,\mu)\to((X\setminus\{e\})_+,\tau_{X\setminus\{e\}})$ and $\eta\colon((X\setminus\{e\})_+,\tau_{X\setminus\{e\}})\to(X,\mu)$ are isomorphisms in $\sQMap$. Indeed, $\zeta\circ\mu=\tau_{X\setminus \{e\}}$ because
\begin{align*}
(\zeta\circ\mu)(\star,\star)&=\bv\limits_{x\in X}\zeta(x,\star)\with\mu(\star,x)\\
&=\bv\limits_{x\in X}\mu^*(x,\star) \with\mu(\star,x)\\
&=(\mu^*\circ\mu)(\star,\star)\\
&=k&\text{(Lemma \ref{zeta-zeta*-k})}\\
&=\tau_{X\setminus \{e\}}(\star,\star)
\end{align*}
and
\begin{align*}
(\zeta\circ\mu)(\star,a)&=\bv\limits_{x\in X}\zeta(x,a)\with\mu(\star,x)\\
&=\bv\limits_{x\in X}\mu(\star,x)\with\mu^*(x\cdot a,\star)\\
&=\bot&\text{(\ref{group-e} and Lemma \ref{zeta-zeta*-k})}\\
&=\tau_{X\setminus \{e\}}(\star,a)
\end{align*}
for all $a\in X\setminus\{e\}$. Meanwhile, $\eta\circ\tau_{X\setminus\{e\}}=\mu$ because
\[(\eta\circ\tau_{X\setminus \{e\}})(\star,x)=\bv\limits_{a\in (X\setminus \{e\})_+}\eta(a,x)\with\tau_{X\setminus \{e\}}(\star,a)=\eta(\star,x)\with k=\eta(\star,x)=\mu(\star,x)\]
for all $x\in X$, where the second equality follows from \eqref{tau_X-def}.
\end{proof}




Lemma \ref{K-sur} guarantees that $\sK\colon\QParMap\to\sQMap$ is an equivalence of categories:

\begin{thm} \label{QMap-iso-sQMap}
Assuming the axiom of choice, $\QParMap$ is equivalent to $\sQMap$ and, therefore, is monadic over $\QMap$.
\end{thm}

As pointed out by an anonymous referee, it is worth exploring whether the conclusion of Theorem \ref{QMap-iso-sQMap} (or Lemma \ref{K-sur}) implies the axiom of choice. However, we are unable to determine whether it is true. So we leave it as a question:

\begin{ques}
Is the equivalence between $\QParMap$ and $\sQMap$ necessarily implies the axiom of choice?
\end{ques}

\section*{Acknowledgement}

The authors acknowledge the support of National Natural Science Foundation of China (No. 12071319) and the Fundamental Research Funds for the Central Universities (No. 2021SCUNL202). We are grateful for helpful remarks received from an anonymous referee which improve the presentation of this paper.






\end{document}